\documentclass[11pt, reqno]{amsart}

\usepackage[OT2, T1]{fontenc}
\usepackage{url}
\usepackage{amsmath}
\usepackage{array}
\usepackage{graphicx}
\usepackage{amsfonts}
\usepackage{amssymb}
\usepackage{amsthm}
\usepackage{hyperref}
\usepackage{colonequals}	
\usepackage{color}
\usepackage{mathabx}		
\usepackage{mathrsfs}
\usepackage{amscd}
\usepackage[all,cmtip]{xy}	
\usepackage{sseq}			
\usepackage{verbatim}
\usepackage{parskip}
\usepackage[in]{fullpage}
\usepackage{mathrsfs} 			
\usepackage{bm} 				
\usepackage{cleveref}
\usepackage{enumitem}
\usepackage{moreenum}
\usepackage[alphabetic,lite]{amsrefs} 	
\usepackage{stmaryrd}	

\DeclareSymbolFont{cyrletters}{OT2}{wncyr}{m}{n}
\DeclareMathSymbol{\Sha}{\mathalpha}{cyrletters}{"58}

\usepackage{color}

\newcommand{\bF}{\mathbb{F}}
\newcommand{\bG}{\mathbb{G}}

\newcommand{\bQ}{\mathbb{Q}}
\newcommand{\bR}{\mathbb{R}}

\newcommand{\bZ}{\mathbb{Z}}

\newcommand{\cA}{\mathcal{A}}

\newcommand{\cF}{\mathcal{F}}
\newcommand{\cG}{\mathcal{G}}

\newcommand{\cN}{\mathcal{N}}
\newcommand{\cO}{\mathcal{O}}

\newcommand{\ra}{\rightarrow}

\newcommand{\xra}{\xrightarrow}

\newcommand{\hra}{\hookrightarrow}

\newcommand{\I}{^{\infty}}

\newcommand{\wh}{\widehat}

\newcommand{\pr}{^{\prime}}
\newcommand{\ce}{\colonequals}
\newcommand{\ov}{\overline}

\renewcommand{\b}{\textbf}

\newcommand{\tensor}{\otimes} 		

\newcommand{\fppf}{\mathrm{fppf}}		
\newcommand{\et}{\mathrm{\acute{e}t}}	

\renewcommand{\i}{^{-1}}

\renewcommand{\implies}{\Rightarrow}

\providecommand{\p}[1]{\left(#1\right)}

\providecommand{\f}[2]{\frac{#1}{#2}}
\DeclareMathOperator{\Ker}{Ker}			
\DeclareMathOperator{\Coker}{Coker}		
\DeclareMathOperator{\im}{Im}			
\DeclareMathOperator{\Spec}{Spec}		
\DeclareMathOperator{\Hom}{Hom}			
\DeclareMathOperator{\Char}{char}		
\DeclareMathOperator{\Br}{Br}		
\DeclareMathOperator{\Gal}{Gal}	
\DeclareMathOperator{\inv}{inv}	
\DeclareMathOperator{\ord}{ord}	
\DeclareMathOperator{\rk}{rk}		
\DeclareMathOperator{\Sel}{Sel}		
\DeclareMathOperator{\Pic}{Pic}		
\DeclareMathOperator{\Frob}{Frob}		
\newcommand{\ba}{\begin{aligned}}
\newcommand{\ea}{\end{aligned}}
\newcommand{\be}{\begin{equation}}
\newcommand{\ee}{\end{equation}}
\newcommand{\pf}{\begin{proof}}
\newcommand{\bpf}{\begin{proof}}
\newcommand{\epf}{\end{proof}}
\newcommand{\bthm}{\begin{thm}}
\newcommand{\ethm}{\end{thm}}
\newcommand{\bprop}{\begin{prop}}
\newcommand{\eprop}{\end{prop}}
\newcommand{\bcor}{\begin{cor}}
\newcommand{\ecor}{\end{cor}}
\newcommand{\brem}{\begin{rem}}
\newcommand{\erem}{\end{rem}}
\newcommand{\brems}{\begin{rems} \hfill \begin{enumerate}[label=\b{\thesubsection.},ref=\thesubsection]}
\newcommand{\remi}{\addtocounter{subsection}{1} \item}
\newcommand{\erems}{\end{enumerate} \end{rems}}
\newcommand{\blem}{\begin{lemma}}
\newcommand{\elem}{\end{lemma}}
\newcommand{\bconj}{\begin{conj}}
\newcommand{\econj}{\end{conj}}
\newcommand{\benum}{\begin{enumerate}[label={(\alph*)}]}
\newcommand{\eenum}{\end{enumerate}}
\newcommand{\bc}{}
\newcommand{\beg}{\begin{eg}}
\newcommand{\eeg}{\end{eg}}
\newcommand{\bcl}{\begin{claim}}
\newcommand{\ecl}{\end{claim}}
\newcommand{\lab}{\label}

\theoremstyle{plain}
\newtheorem{thm}[subsection]{Theorem}
\Crefname{thm}{Theorem}{Theorems}

\Crefname{rethm}{Theorem}{Theorem}
\newtheorem{prop}[subsection]{Proposition}
\Crefname{prop}{Proposition}{Propositions}

\Crefname{q}{Question}{Questions}
\Crefname{eg}{Example}{Examples}

\Crefname{Problem}{Problem}{Problems}
\newtheorem{conj}[subsection]{Conjecture}
\Crefname{conj}{Conjecture}{Conjectures}
\newtheorem{cor}[subsection]{Corollary}
\Crefname{cor}{Corollary}{Corollaries}

\newtheorem{lemma}[subsection]{Lemma}

\theoremstyle{remark}
\newtheorem{claim}[equation]{Claim}
\Crefname{claim}{Claim}{Claims}

\theoremstyle{definition}
\newtheorem{eg}[subsection]{Example}
\newtheorem{rem}[subsection]{Remark}
\Crefname{rem}{Remark}{Remarks}
\newtheorem*{rems}{Remarks}

\newtheoremstyle{subsection-tweak}
   {11pt}
   {3pt}%
   {}
   {}%
   {\bfseries}
   {}%
   {.5em}
   {\thmnumber{\@{#1}{}\@{#2}.}%
    \thmnote{~{\bfseries#3.}}}

\Crefname{innercustomconj}{Conjecture}{Conjecture}

\theoremstyle{subsection-tweak}
\newtheorem{pp}[subsection]{}
\newcommand{\bpp}{\begin{pp}}
\newcommand{\epp}{\end{pp}}

\numberwithin{equation}{subsection}

\begin{document}
\author{K\k{e}stutis \v{C}esnavi\v{c}ius}
\title{Selmer groups and class groups}
\date{\today}
\subjclass[2010]{Primary 11G10; Secondary 11R23, 11R29, 11R58}
\keywords{Selmer group, class group, fppf cohomology, Iwasawa theory}
\address{Department of Mathematics, Massachusetts Institute of Technology, Cambridge, MA 02139, USA}
\email{kestutis@math.mit.edu}
\urladdr{http://math.mit.edu/~kestutis/}

\begin{abstract}  
Let $A$ be an abelian variety over a global field $K$ of characteristic $p \ge 0$. If $A$ has nontrivial (resp.~full) $K$-rational $l$-torsion for a prime $l \neq p$, we exploit the fppf cohomological interpretation of the $l$-Selmer group $\Sel_l A$ to bound $\#\Sel_l A$ from below (resp.~above) in terms of the cardinality of the $l$-torsion subgroup of the ideal class group of $K$. Applied over families of finite extensions of $K$, the bounds relate the growth of Selmer groups and class groups. For function fields, this technique proves the unboundedness of $l$-ranks of class groups of quadratic extensions of every $K$ containing a fixed finite field $\bF_{p^n}$ (depending on $l$). For number fields, it suggests a new approach to the Iwasawa $\mu = 0$ conjecture through inequalities, valid when $A(K)[l] \neq 0$, between Iwasawa invariants governing the growth of Selmer groups and class groups in a $\bZ_l$-extension.
\end{abstract}

\maketitle



\section{Introduction}

Fix a prime $l$, a number field $K$, an abelian variety $A \ra \Spec K$ of dimension $g > 0$, and let $L/K$ range in some family of finite extensions. Our goal is to relate, in favorable situations, the growth of the $l$-torsion subgroup $\Pic(\cO_L)[l]$ of the ideal class group of $L$ and that of the $l$-Selmer group $\Sel_l A_L$. Concrete expectations in the case of quadratic $L/K$ are provided by folklore conjectures:

\bconj \lab{conj-cl}
As $L/K$ ranges over quadratic extensions, $\#\Pic(\cO_L)[l]$ is unbounded.
\econj

\bconj \lab{conj-sel}
As $L/K$ ranges over quadratic extensions, $\#\Sel_l A_L$ is unbounded.
\econj

\brems 
\remi\lab{c-c-kn}
\Cref{conj-cl} is known for $l = 2$ due to the genus theory of Gauss, but is open for every pair $(K, l)$ with $l$ odd; in the $K = \bQ$ case, much more precise predictions are available through the Cohen--Lenstra heuristics \cite{CL84}. The conjectured (but not universally believed) unboundedness of $\rk A(L)$ would imply \Cref{conj-sel}, which is known for $l = 2$ if $g = 1$ \cite{CS10}*{Thm.~3}\footnote{The case when $A$ does not have potential complex multiplication is due to B\"{o}lling \cite{Bol75}*{pp.~170--171}. Both papers concern the (stronger) unboundedness of cardinalities of $2$-torsion subgroups of Shafarevich--Tate groups.} and for $l = 2$ in certain $g > 1$ cases (see Remarks \ref{our-contr} and \ref{lit-rem}), but is open for every pair $(A, l)$ with $l$ odd. 

\remi 
If \Cref{conj-cl} (resp.,~\ref{conj-sel}) is known for $(K, l)$, it follows for $(K\pr, l)$ for every finite extension $K\pr/K$, see \Cref{pic-extn} (resp.,~\ref{sel-extn}).
\erems

We relate the conjectures by proving their equivalence after replacing $K$ by a finite extension:

\bthm[\Cref{sel-pic}] \lab{quad-eq} \hfill
\benum
\item \lab{quad-eq-a}
If $A$ has $\bZ/l\bZ$ or $\mu_l$ as a $K$-subgroup, then \Cref{conj-cl} for $K$ implies \Cref{conj-sel} for~$A$.

\item \lab{quad-eq-b}
If $A[l]$ has a filtration by $K$-subgroups with subquotients isomorphic to $\bZ/l\bZ$ or $\mu_l$, then \Cref{conj-sel} for $A$ implies \Cref{conj-cl} for $K$.
\eenum
\ethm

\brems
\remi \lab{our-contr}
The known $l = 2$ case of \Cref{conj-cl} therefore proves the $l = 2$ and $A(K)[2] \neq 0$ case of \Cref{conj-sel}. Restricting further to $g = 1$, this combines with the unboundedness of $\#\Sel_2 A_L$ proved by Klagsbrun, Mazur, and Rubin \cite{Kla11}*{1.2} 
under the $A(K)[2] = 0$ assumption to reprove \Cref{conj-sel} in the $(g, l) = (1, 2)$ case.

\remi
Even though the idea that Selmer groups and class groups are related is not new (compare, e.g.,~\cite{Sch96}), the relationship furnished by \Cref{quad-eq} is sharper than those available previously. Moreover, it is specific neither to quadratic $L/K$ nor to number fields: \S\ref{bd},~containing its proof, works in the setting of bounded degree extensions of any fixed global field $K$.
\erems

\bpp[The method of the proof] \lab{outline}
Under the assumptions of \ref{quad-eq-a} (resp.,~\ref{quad-eq-b}) of \Cref{quad-eq}, we prove lower (resp.,~upper) bounds for $\#\Sel_l A$ in terms of $\#\Pic(\cO_K)[l]$ in \S\ref{low} (resp.,~\S\ref{up}), which we apply after base change to $L$. As for the bounds themselves, the fppf cohomological interpretation of Selmer groups provides the idea. To explain it, assume for simplicity that $A[l] \cong (\bZ/l\bZ)^g \oplus \mu_l^g$ over $K$, and let $S$ be the spectrum of the ring of integers of $K$ and $\cA \ra S$ the N\'{e}ron model of $A$. The N\'{e}ron property of $\cA[l]_{S[\f{1}{l}]}$ \cite{Ces13}*{B.6} forces $\cA[l]_{S[\f{1}{l}]} \cong (\bZ/l\bZ)^g \oplus \mu_l^g$. Passing to cohomology, both $\#H^1(S[\f{1}{l}], \bZ/l\bZ)$ and $\#H^1(S[\f{1}{l}], \mu_l)$ relate to $\#\Pic(S)[l]$ (see \Cref{const,mu}), whereas $H^1(S[\f{1}{l}], \cA[l]) \subset H^1(K, A[l])$ is defined by local conditions \cite{Ces13}*{4.2}, which at finite places of good reduction agree with those defining $\Sel_l A \subset H^1(K, A[l])$ \cite{Ces13}*{2.5}; it remains to quantify the resulting relation between $\#H^1(S[\f{1}{l}], \cA[l])$ and $\#\Sel_l A$. 
\epp

\bpp[The function field case] \lab{ff-ano}
The argument sketched in \ref{outline} continues to work for a global function field $K$ of positive characteristic $p \neq l$. For such $K$, the analogue of \Cref{conj-sel} is known in the case of a constant supersingular elliptic curve: $\rk A(L)$ is unbounded due to the work of Shafarevich and Tate \cite{TS67}. With this input, we prove the analogue of \Cref{conj-cl} for every $K$ containing a fixed finite field $\bF_{p^n}$ (depending on $l$) and consequently, for such $K$, also the analogue of \Cref{conj-sel} for $A$ that have $\bZ/l\bZ$ or $\mu_l$ as a $K$-subgroup. For precise statements, see \Cref{S-T} and \Cref{sel-unb-ff}. As in the number field case discussed in Remark \ref{c-c-kn}, no case of the analogue of \Cref{conj-cl} was previously known for odd $l$ (for $l = 2$, see \cite{Mad72}*{Thm.~3}).
\epp

\bpp[Applications to Iwasawa theory]
The bounds mentioned in \ref{outline} lead to inequalities of \Cref{Zp-up,Zp-low} between the Iwasawa invariants governing the growth of Selmer groups and class groups in the layers of a $\bZ_p$-extension. These inequalities imply our main result concerning Iwasawa theory (for a detailed discussion and other results see \S\S\ref{Iw-prel}--\ref{cycl-concl}):
\epp

\bthm[\Cref{main-Iw}] \lab{main-Iw-ano}
For a prime $p$ and a number field $K$, to prove the Iwasawa $\mu = 0$ conjecture for the cyclotomic $\bZ_p$-extension $K_\infty/K$, it suffices to find an abelian $K$-variety $A$ such~that
\begin{enumerate}[label={(\roman*)}]
\item \lab{main-Iw-i}
$A$ has good ordinary reduction at all places above $p$,

\item
$A$ has $\bZ/p\bZ$ as a $K$-subgroup, 

\item \lab{main-Iw-iii}
$\Hom(\Sel_{p\I} A_{K_\infty}, \bQ_p/\bZ_p)$ is a torsion module over the Iwasawa algebra and has $\mu$-invariant~$0$.
\eenum
\ethm

\brem
In fact, it suffices to find such an $A$ after replacing $K$ by a finite extension, see \Cref{mu-pic-gr}. It is not clear, however, how to take advantage of the apparent flexibility of choice: for arbitrary $K$ and $p$, \ref{main-Iw-iii} alone seems nontrivial to fulfill. For $K = \bQ$ and $p = 5$, the elliptic curve 11A3 satisfies \ref{main-Iw-i}--\ref{main-Iw-iii} \cite{Gre99}*{pp.~120--124}; with this $A$, \Cref{main-Iw-ano} reproves an easy case of the Ferrero--Washington theorem (which is not used in loc.~cit., so the argument is not circular).
\erem

\bpp[The contents of the paper]
The bounds discussed in \ref{outline} are essential for all subsequent applications and are proved in \S\S\ref{low}--\ref{up}. These technical sections rely on (standard but crucial) auxiliary computations of appendices \ref{app-a} and \ref{app-b}. \Cref{quad-eq} is proved in \S\ref{bd}, which applies the inequalities of \S\S\ref{low}--\ref{up} in families of bounded degree extensions of $K$. Both \S\S\ref{low}--\ref{bd} and the appendix \ref{app-b} work under the assumption that $K$ is a global field. Special cases of function field analogues of \Cref{conj-cl,conj-sel} are proved in \S\ref{evidence}. The remaining \S\S\ref{Iw-prel}--\ref{cycl-concl} discuss Iwasawa theory (and assume that $K$ is a number field). The introductory \S\ref{Iw-prel} records how Iwasawa invariants control the growth of $\Pic(\cO_K)[p^m]$ and $\Sel_{p^m} A$; this deviates from the standard discussion that concerns $\Pic(\cO_K)[p\I]$ and $\Sel_{p\I} A$. Inequalities between Iwasawa invariants of class groups and Selmer groups result from the bounds of \S\S\ref{low}--\ref{up} and are the subject of \S\ref{Iw}. The final \S\ref{cycl-concl} summarizes the conclusions for the cyclotomic $\bZ_p$-extension (\S\S\ref{Iw-prel}--\ref{Iw} allow an arbitrary $\bZ_p$-extension).
\epp

\bpp[Notation] \lab{not} 
The notation set in this paragraph is in place for the rest of the paper; deviations, if any, are recorded in the beginning of each section. Let $l$ be a prime, $m$ a positive integer, and $K$ a global field. If $\Char K = 0$, let $S$ be the spectrum of the ring of integers of $K$; if $\Char K > 0$, let $S$ be the smooth proper curve over a finite field such that the function field of $S$ is $K$. Let $v$ be a place of $K$ and $K_v$ the corresponding completion; if $v\nmid \infty$, then $v$ identifies with a closed point of $S$, and $\cO_v$ and $\bF_v$ denote the ring of integers and the residue field of $K_v$. Let $r_1$ and $r_2$ be the number of real and complex places of $K$. Let $A \ra \Spec K$ be an abelian variety of dimension $g > 0$ and $\cA \ra S$ its N\'{e}ron model. For $v\in S$, let $\Phi_v$ be the \'{e}tale $\bF_v$-group scheme of connected components of $\cA_{\bF_v}$. For a finite extension $L/K$, the formation of $S$, $\cA$, $\Phi_v$ is \emph{not} compatible with base change, and we denote by $S^L$, $\cA^L$, $\Phi_w^L$ their analogues over $L$ (note that $S^L$ is the normalization of $S$ in $L$). 
\epp

\bpp[Conventions] \lab{conv}
To simplify the computations, $\gtrsim_{A, L, \dotsc}$ and $\lesssim_{A, L, \dotsc}$ denote inequalities up to implied constants that depend only on the indicated parameters (note that $A$, being a morphism $A \ra \Spec K$, includes dependence on $K$); when no parameters are indicated, the ones used last are taken. Also, $\sim$ stands for ``$\gtrsim$ and $\lesssim$''. When needed (e.g., for forming composita or intersections), a choice of a separable closure $\ov{F}$ of a field $F$ is made implicitly (and compatibly for overfields). The \'{e}tale fundamental group of an integral scheme is based at a geometric generic point. Fppf cohomology is denoted by $H^i$; when the coefficient sheaf is a smooth group scheme, the identification with \'{e}tale cohomology \cite{Gro68}*{11.7 1$^{\circ}$)} is implicit and similarly for further identifications with Galois cohomology. Fppf cohomology with compact supports that takes into account infinite primes \cite{Mil06}*{III.0.6 (a)} is denoted by $H^i_c$. All quotients are taken in the big fppf topos, and $X_\fppf$ denotes the big fppf site of the scheme $X$. The $l^m$-Selmer group $\Sel_{l^m} A$ is the preimage of $\prod_v A(K_v)/l^mA(K_v) \subset \prod_v H^1(K_v, A[l^m])$ in $H^1(K, A[l^m])$, where fppf cohomology is necessary if $l = \Char K$. For a nonempty open $U \subset S$, the number of closed points of $S$ not in $U$ is $\#(S\setminus U)$. If $\Char K = 0$, then $\Pic_+(S)$ is the narrow ideal class group of $K$; if $\Char K > 0$, then $\Pic_+(S) \ce \Pic(S)$. For an integer $n$ and a scheme $X$, the open subscheme on which $n$ is invertible is $X[\f{1}{n}]$.
\epp

\subsection*{Acknowledgements} I thank Bjorn Poonen for many helpful discussions and suggestions. I thank Julio Brau, Pete Clark, Tim Dokchitser, Jordan Ellenberg, Zev Klagsbrun, Barry Mazur, Filip Najman, Karl Rubin, Doug Ulmer, Jeanine van Order, Larry Washington, and David Zureick-Brown for helpful conversations or correspondence regarding the material of the paper. I thank the referee for helpful suggestions. Part of the research presented here was carried out during the author's stay at the Centre Interfacultaire Bernoulli (CIB) in Lausanne during the course of the program ``Rational points and algebraic cycles''. I thank CIB, NSF, and the organizers of the program for a lively semester and the opportunity to take part.

\section{Lower bounds for Selmer groups in terms of class groups} \lab{low}

Mimicking \cite{Mil06}*{before II.3.4}, for a nonempty open $U \subset S$ and a sheaf $\cF$ on $U_{\fppf}$, we define 
\[
D^1(U, \cF) \ce \im(H^1_c(U, \cF) \ra H^1(U, \cF)).
\]

\bprop \lab{D1-Sel}
If $U \subset S$ is a nonempty open subscheme for which $A$ has semiabelian reduction at all $v\in U$ with $\Char \bF_v = l$, then 
\[
\xymatrix{
  D^1(U, \cA[l^m]) \ar@{^(->}[r] \ar[d] & H^1(K, A[l^m]) \ar[d] \\
  \prod_{v \in U} H^1(\cO_v, \cA[l^m]) \times \prod_{\substack{v\not\in U}} 0 \ar@{^(->}[r] &\prod_v  H^1(K_v, A[l^m]),
}
\]
is Cartesian. If, moreover, $l \neq \Char K$ or $U = S$, then, taking intersections inside $H^1(K, A[l^m])$,
\[ \ba
\#\p{\f{D^1(U, \cA[l^m])}{D^1(U, \cA[l^m]) \cap \Sel_{l^m} A}} &\le \prod_{v\in U} \f{\#\Phi_v(\bF_v)}{\#(l^m\Phi_v)(\bF_v)},\\
\#\p{\f{\Sel_{l^m} A}{D^1(U, \cA[l^m]) \cap \Sel_{l^m} A}} &\le \prod_{v\in U} \f{\#\Phi_v(\bF_v)}{\#(l^m\Phi_v)(\bF_v)} \cdot \prod_{v\in S\setminus U} \p{l^{mg[K_v : \bQ_l]}\cdot \#A(K_v)[l^m]} \cdot \prod_{\substack{\text{real }v \\ l = 2}} \# \pi_0(A(K_v)),
\ea\]
where $[K_v : \bQ_l] \ce 0$ unless $K_v$ is a finite extension of $\bQ_l$.
\eprop

\bpf
For the diagram, use the similar description of $H^1(U, \cA[l^m]) \subset H^1(K, A[l^m])$ \cite{Ces13}*{4.2~and~B.5} and the compactly supported cohomology exact sequence 
\cite{Mil06}*{III.0.6~(a)}. For the inequalities, compare the defining local conditions by means of \cite{Ces13}*{2.5 (a)} and \Cref{card-kum}.
\epf

\bthm  \lab{low-bds}
Suppose that $A[l^m]$ has a $K$-subgroup $G \cong \bigoplus_{i\in I} \bZ/l^{a_i}\bZ \oplus \bigoplus_{j \in J} \mu_{l^{b_j}}$ with $a_i, b_j \ge 1$.

\benum
\item \lab{low-bds-a}
Set $r \ce r_1$ if $l = 2$, and $r \ce 0$ if $l \neq 2$; also $\{v\mid l \} \ce \emptyset$ if $\Char K > 0$. If $l \neq \Char K$, then
\[
\#\Sel_{l^m} A \gtrsim_{g, l, m}  \f{\prod_i  \#\Pic(S[\tfrac{1}{l}])[l^{a_i}] \prod_{j} \#\Pic_+(S)[l^{b_j}]}{2^{r\cdot \#J}\cdot l^{r_2\sum_j b_j} \cdot \prod_{v \in S[\f{1}{l}]} \f{\#\Phi_v(\bF_v)}{\#(l^m\Phi_v)(\bF_v)}\cdot \prod_j \prod_{v\mid l} \#\mu_{l^{b_j}}(K_v)}.
\]

\item \lab{low-bds-b}
If $J = \emptyset$ and $A$ has semiabelian reduction at all $v$ with $\Char \bF_v = l$, then
\[\ba 
\#\Sel_{l^m} A \gtrsim_{g, l, m} \f{\prod_i \#\Pic(S)[l^{a_i}]}{  \prod_{v\nmid \infty} \f{\#\Phi_v(\bF_v)}{\#(l^m\Phi_v)(\bF_v)}} .
\ea \]
\eenum 
\ethm

\bpf 
We give the similar proofs together. For \ref{low-bds-a}, set $U \ce S[\f{1}{l}]$; for \ref{low-bds-b}, set $U \ce S$. By \Cref{D1-Sel}, $\#\Sel_{l^m} A \ge \#D^1(U, \cA[l^m]) \cdot \p{ \prod_{v\in U} \f{\#\Phi_v(\bF_v)}{\#(l^m\Phi_v)(\bF_v)} }\i$. Let $\cG \ra U$ be the group smoothening of the schematic image of $G \ra \cA_{U}$; by \cite{BLR90}*{7.1/6}, $\cG$ is the N\'{e}ron model of $G$, hence $\cG \cong \bigoplus_i \bZ/l^{a_i}\bZ \oplus \bigoplus_j \mu_{l^{b_j}}$. The $U$-homomorphism $\cG \xra{f} \cA[l^m]$ has generic fiber $G \hra A[l^m]$; moreover, $H^1(U, \cG) \subset H^1(K, G)$ and $H^1(U, \cA[l^m]) \subset H^1(K, A[l^m])$ \cite{Ces13}*{A.5~and~B.5}. Therefore, $\#\Ker H^1(f) \lesssim 1$, giving $\#D^1(U, \cA[l^m]) \gtrsim \#D^1(U, \cG)$. The conclusion follows by combining the obtained inequalities with \Cref{const-c,mu-c} and the exact sequence \cite{Mil06}*{III.0.6 (a)}.
\epf

\section{Upper bounds for Selmer groups in terms of class groups} \lab{up}

Assume in this section that $l \neq \Char K$. Contrary to the lower bounds in \Cref{low-bds}, we do not use implied constants in the upper bounds in \Cref{up-bds}. This makes the inequalities less pleasant but has the advantage of providing explicit lower bounds on the cardinalities of $l$-torsion subgroups of class groups when \Cref{up-bds} is applied to an abelian variety of high rank. For instance, one may hope for a practical approach to \Cref{quad-eq}: by finding an elliptic curve $E \ra \Spec \bQ$ for which $E(\bQ)[l] \neq 0$ with $l$ odd and a quadratic $F/\bQ$ for which $\rk E(F)$ is large, one would get a quadratic number field with large class group $l$-rank $r_l \ce \dim_{\bF_l} \Pic(S^F)[l]$. The current records (among quadratic $F$) $r_3 = 6$ \cite{Que87} and $r_5 = 4$ \cite{Sch83} exploit relations with elliptic curves.

\bthm \lab{up-bds}
Fix a nonempty open $U \subsetneq S[\f{1}{l}]$ for which $\cA_U \ra U$ is an abelian scheme. Set $r \ce r_1$ if $l = 2$, and $r \ce 0$ if $l \neq 2$; also $[K : \bQ] \ce 0$ if $\Char K > 0$. If $A[l^m]$ has a filtration by $K$-subgroups $N_j$ with subquotients isomorphic to $\bZ/l^{a_i}\bZ$ or $\mu_{l^{b_j}}$ with $a_i, b_j \ge 1$,~then
\[ \ba
\#\Sel_{l^m} A \le &\prod_i \#(\Pic_+ S / l^{a_i} \Pic_+ S) \prod_j \#\Pic(U)[l^{b_j}] \cdot l^{[K: \bQ]\sum_i a_i + (r_1 + r_2 + \#(S\setminus U) - 1)\sum_j b_j} \cdot \\ &\prod_j \#\mu_{l^{b_j}}(K) \cdot \prod_i \prod_{v\in S \setminus U} \#\mu_{l^{a_i}}(K_v),
\ea \] 
and also
\[ \ba
\# \Sel_{l^m} A \le &\prod_{i} \#(\Pic U/l^{a_i} \Pic U)  \prod_j \#(\Pic_+ S / l^{b_j} \Pic_+ S) \cdot l^{mg[K: \bQ] + (\#(S \setminus U) - 1)\sum_i a_i + (r_1 + r_2 - 1)\sum_j b_j} \cdot \\ & \prod_i 2^{r} \cdot \prod_{\substack{\text{real }v \\ l = 2}} \# \pi_0(A(K_v))\i \cdot \#A[l^m](K) \cdot  \prod_j \prod_{v\in S \setminus U} \#\mu_{l^{b_j}}(K_v).
\ea \]
\ethm

\bpf 
Let $\cN_j$ be the schematic image of $N_j \ra \cA[l^m]_{U}$. By \cite{EGAI}*{9.5.5--6}, \cite{EGAIV2}*{2.8.5--6}, \cite{TO70}*{p.~17 Lemma 5}, and finiteness of $\cA[l^m]_U$, the $\cN_j$ filter $\cA[l^m]_U$ by finite \'{e}tale $U$-subgroups. Due to finiteness, the \'{e}tale subquotients $\cN_{j + 1}/\cN_j$ are the N\'{e}ron models of the $N_{j + 1}/N_j$ and hence identify with $\bZ/l^{a_i}\bZ$ or $ \mu_{l^{b_j}}$. Therefore, \Cref{const,mu,const-c,mu-c} bound $\#H^1(U, \cA[l^m])$ and $\#H^1_c(U, \cA[l^m])$ through cohomology sequences, and the claimed inequalities follow by combining these bounds with the following observations:
\begin{enumerate}[label={(\roman*)}]
\item 
For the first inequality: by \cite{Ces13}*{2.5 (d) and 4.2}, $\#\Sel_{l^m} A \le \#H^1(U, \cA[l^m])$;

\item 
For the second: by \cite{Mil06}*{III.0.6 (a)} and \Cref{D1-Sel}, writing $\wh{H}^i$ for Tate cohomology,
\[
\#D^1(U, \cA[l^m]) \le \#H^1_c(U, \cA[l^m]) \cdot \prod_{v \in S\setminus U} \#A(K_v)[l^m]\i \cdot \prod_{v \mid \infty} \#\wh{H}^0(K_v, A[l^m])\i \cdot \#A[l^m](K),\text{ and}
\]
\[
\f{\#\Sel_{l^m} A}{\#D^1(U, \cA[l^m])} \le l^{mg[K : \bQ]} \cdot \prod_{v\in S\setminus U} \#A(K_v)[l^m] \cdot \prod_{\substack{\text{real }v \\ l = 2}} \# \pi_0(A(K_v));
\]
moreover, if $l = 2$ and $v$ is real, by \Cref{card-kum}\ref{card-kum-c} and \cite{GH81}*{1.3},
\[
\#\wh{H}^0(K_v, A[l^m]) = \#H^1(K_v, A[l^m]) = \#\pi_0(A(K_v))^2. \qedhere
\]
\eenum
\epf

\brems
\remi \lab{ch-a}
The two bounds are incomparable in general; they yield different bounds in \Cref{Zp-up}. 

\remi \lab{ch-b}
When $\bZ/l^{a_i}\bZ \cong \mu_{l^{a_i}}$ over $K$, the two interpretations of the corresponding subquotient result in different right hand sides of the inequalities of \Cref{up-bds}, and hence also in the flexibility of choosing the best bound. Similarly for $\mu_{l^{b_j}}$.
\erems

\section{Growth of Selmer groups and class groups in extensions of bounded degree} \lab{bd}

\bthm \lab{bdd-deg}
Let $L/K$ be an extension of degree at most $d$.
\benum
\item \lab{bdd-deg-a}
If either
\begin{enumerate}[label={(\roman*)}] 
\item \lab{a1}
$A$ has $\bZ/l\bZ$ or $\mu_l$ as a $K$-subgroup, and $l \neq \Char K$, or
\item \lab{a2}
$A$ has everywhere semiabelian reduction and $\bZ/l\bZ$ as a $K$-subgroup,
\eenum
then
\[
\#\Sel_{l^m} A_L \gtrsim_{A, d, l} \#\Pic(S^L)[l].
\]

\item \lab{bdd-deg-b}
If $l \neq \Char K$ and $A[l]$ has a filtration with subquotients isomorphic to $\bZ/l\bZ$ or $ \mu_{l}$, then
\[
\#\Sel_{l} A_L \lesssim_{A, d, l} \#\Pic(S^L)[l]^{2g}.
\]
\eenum
\ethm

\bpf \hfill
\benum
\item 
This follows from \Cref{low-bds} since, letting $w$ denote a place of $L$, we have
\begin{enumerate}[label={(\arabic*)}]
\item \lab{triv-invp}
$\#\Pic (S^L[\f{1}{l}])[l] \sim_{K, d, l} \#\Pic (S^L)[l]$ if $\Char K \neq l$, because $\#(S^L\setminus S^L[\f{1}{l}])$ is bounded;
\item \lab{triv-real}
$\#\Pic_+ (S^L)[l] \sim_{K, d} \#\Pic (S^L)[l]$, because the number of real $w$ is bounded;
\item \lab{triv-bad}
There is a bounded number of $w$'s of bad reduction for $A$; moreover, for each such $w$, 
\begin{enumerate}[label={(\greek*)}]
\item 
If $\Char K = 0$, up to isomorphism there are only finitely many possibilities for~$A_{L_w}$.

\item
In general, $\f{\#\Phi_w^L(\bF_w)}{\#(l\Phi_w^L)(\bF_w)} \le \# \Phi_w^L[l]$, and, if $l \neq \Char \bF_w$ or the reduction is semiabelian, then $\# \Phi_w^L[l] \le \#\cA^L[l]_{\bF_w} \sim_{g, l} 1$, as is seen by inspecting the finite part \cite{EGAIV4}*{18.5.11 c)} of the quasi-finite separated $\cA^L[l]_{\cO_w}$.
\eenum
\eenum
\item 
This follows from (either part of) \Cref{up-bds}: one argues as in \ref{triv-invp} and \ref{triv-real} and uses
\begin{enumerate}[label={(\arabic*)}] \addtocounter{enumii}{3}
\item 
$\#(\Pic_+ S^L/l \Pic_+ S^L) \sim_{K, d, l} \#\Pic (S^L)[l]$. \qedhere
\eenum
\eenum
\epf

\bcor \lab{cycl-p}
If either \ref{a1} or \ref{a2} of \Cref{bdd-deg}\ref{bdd-deg-a} hold, then $\#\Sel_l A_{L}$ is unbounded as $L/K$ ranges over degree $l$ extensions.
\ecor

\bpf
Indeed, $\#\Pic (S^L)[l]$ is unbounded \cite{Mad72}*{Thm.~3}.
\epf

\bcor
If $l \neq \Char K$, then $\#\Sel_l A_L$ is unbounded as $L/K$ ranges over extensions of degree at most $l^{2g + 1} - l$.
\ecor

\bpf
Indeed, $A$ acquires a nontrivial $l$-torsion point over an extension of degree at most $l^{2g} - 1$.
\epf

\brem \lab{lit-rem}
There are several results in the literature concerned with proving the unboundedness of $\#\Sha(A_L)[l]$ (and hence that of $\#\Sel_l A_L$) as $L$ ranges over degree $l$ extensions of $K$: \cite{CS10}*{Thm.~3} treats the case $\dim A = 1$ and $l \neq \Char K$, whereas \cite{Cre11}*{Thm.~1.1}, improving \cite{Cla04}*{Thm.~7}, allows arbitrary dimension but imposes restrictions (which are satisfied after passing to a finite extension) on the N\'{e}ron--Severi group of $A$. In contrast, \Cref{cycl-p} has no dimension or N\'{e}ron--Severi assumptions but constrains $A[l]$ and only gives Selmer growth.
\erem

If $l \neq \Char K$, the assumptions of \ref{bdd-deg-a} and \ref{bdd-deg-b} in \Cref{bdd-deg} are satisfied after passing to a suitable finite extension $K\pr/K$; standard lemmas \ref{pic-extn} and \ref{sel-extn}, which are also used in \S\ref{Iw}, clarify in \Cref{sel-pic} how this affects the unboundedness questions.

\blem \lab{pic-extn}
Let $L$ be a global field and $L\pr/L$ an extension of degree at most $d$. Then 
\[
\#\Pic (S^{L\pr})[n] \gtrsim_{d, n} \#\Pic (S^L)[n].
\]
\elem

\bpf
For number fields, the claim is clear from the theory of the Hilbert class field: if $H/L$ is an unramified abelian extension with Galois group killed by $n$, then so is $HL\pr/L\pr$, for which $[HL\pr : L\pr] \ge \f{1}{d} [H : L]$. The proof in the function field case is the same---the link to unramified abelian extensions is provided by \Cref{const}\ref{const-a} applied to the prime factors of $n$: $\#\Pic(S^L)[n] \sim_{n} \#H^1(S^L, \bZ/n\bZ) = \#\Hom(\pi_1^{\et}(S^L), \bZ/n\bZ) = [H_L : L]$ where $H_L/L$ is the maximal (in $\ov{L}$) unramified abelian extension with Galois group killed by $n$, and similarly for $L\pr$.
\epf

\blem \lab{sel-extn}
Let $L$ be a global field, $A$ a $g$-dimensional abelian variety over $L$, and $L\pr/L$ an extension of degree at most $d$. If $\Char L \nmid n$, then 
\[
\#\Sel_n A_{L\pr} \gtrsim_{d, g, n} \# \Sel_n A.
\]
\elem

\bpf
Let ${\rm{R}}_{L\pr/L}$ denote the restriction of scalars. By \cite{CGP10}*{A.5.1--2, A.5.4 (1), A.5.7}, 
\be\ba \lab{sel-ba-ch}
\xymatrix{
0 \ar[r] & A[n] \ar[r]\ar@{^(->}[d]^{a} & A \ar[r]^{n}\ar@{^(->}[d] & A \ar[r]\ar@{^(->}[d] & 0 \\
0 \ar[r] & {\rm{R}}_{L\pr/L}(A[n]_{L\pr}) \ar[r] & {\rm{R}}_{L\pr/L}(A_{L\pr}) \ar[r]^{n} & {\rm{R}}_{L\pr/L}(A_{L\pr}) \ar[r] & 0 }
\ea\ee
is a morphism of short exact (in the big \'{e}tale site of $L$) sequences of smooth $L$-group schemes. Moreover, ${\rm{R}}_{L\pr/L}(A[n]_{L\pr})$ is finite \'{e}tale with $\#{\rm{R}}_{L\pr/L}(A[n]_{L\pr})\sim 1$: for separable $L\pr/L$, this is evident after base change to $L\pr$, \cite{CGP10}*{A.5.13} handles the purely inseparable case, and in general one uses the transitivity of ${\rm{R}}_{L\pr/L}$. Consequently, $\#\Ker H^1_\et(a) \sim 1$, and since $H^i_\et(L, {\rm{R}}_{L\pr/L}(A[n]_{L\pr})) \cong H^i_\et(L\pr, A[n])$ \cite{SGA4.5}*{p.~24 II.3.6}, it remains to see that $H^1_\et(a)$ respects the $n$-Selmer subgroups. This is evident from the compatibility of the formation of \eqref{sel-ba-ch} with any base change and the well known $L\pr \tensor_{L} L_v \cong \prod_{w\mid v} L\pr_w$ \cite{Ser79}*{II.\S3 Thm.~1 (iii)} 
for a place $v$ of~$L$.
\epf

\brem
For separable $L\pr/L$, one reduces to the Galois case and applies the inflation-restriction sequence in Galois cohomology to obtain another proof of \Cref{sel-extn}.
\erem

\bcor\lab{sel-pic} Let $L/K$ range in a family of finite extensions of bounded degree.
\benum
\item \lab{sel-pic-a}
For a finite extension $K\pr/K$ for which either \ref{a1} or \ref{a2} of \Cref{bdd-deg}\ref{bdd-deg-a} hold, if $\#\Pic(S^L)[l]$ is unbounded, then so is $\#\Sel_{l} A_{K\pr L}$.

\item \lab{sel-pic-b}
Assume that $l \neq \Char K$. For a finite extension $K\pr/K$ for which $A[l]_{K\pr}$ has a filtration with subquotients isomorphic to $\bZ/l\bZ$ or $\mu_{l}$, if $\#\Sel_{l} A_L$ is unbounded, then so is $\#\Pic (S^{K\pr L})[l]$.
\eenum
\ecor

\bpf
Combine \Cref{bdd-deg} with \Cref{pic-extn,sel-extn}.
\epf

\section{Special cases of the function field analogues of \Cref{conj-cl,conj-sel}} \lab{evidence}

For this section, fix a prime $p$ and suppose that $\Char K = p$, i.e., $K$ is a finite extension of $\bF_p(t)$. The analogues in question assume that $l \neq p$ and predict that $\#\Pic(S^L)[l]$ and $\#\Sel_l A_L$ should be unbounded as $L$ ranges over quadratic extensions of $K$.
We show that this is indeed the case if one replaces $K$ by a finite extension depending on $l$ (and also on $A$ in the Selmer group case). The key input is the work of Shafarevich and Tate \cite{TS67} on unboundedness of ranks of quadratic twists of a constant supersingular elliptic curve.

\bthm \lab{S-T}
For each prime power $l^m$ with $l \neq p$, there is a $q = p^{n(l, m)}$ such that if $\bF_q \subset K$, then the number of $\bZ/l^m\bZ$-summands of $\Pic(S^L)[l^m]$ is unbounded as $L/K$ ranges over quadratic extensions of the form $L = L\pr K$ for quadratic extensions $L\pr/\bF_p(t)$. In particular, with $n \ce n(l, 1)$, the analogue of \Cref{conj-cl} holds for $l$ and every global field containing $\bF_{p^n}$.
\ethm

\bpf
Take a supersingular elliptic curve $E \ra \Spec \bF_p$ (see \cite{Wat69}*{4.1 (5)} for its existence proved by Deuring). Let $q$ be such that $E_{\bF_q}[l^m] \cong \bZ/l^m\bZ \oplus \mu_{l^m}$, and hence also $E_{S^L}[l^m] \cong \bZ/l^m\bZ \oplus \mu_{l^m}$ for every $L$. By \cite{Ces13}*{5.4 (c)}, $H^1(S^L, E_{S^L}[l^m]) = \Sel_{l^m} E_L$, and by the result of Shafarevich and Tate \cite{Ulm07}*{1.4}, $\rk E(L)$ and hence also the number of $\bZ/l^m\bZ$-summands of $\Sel_{l^m} E_L$ are unbounded. It remains to note that by the proofs of \Cref{const,mu}, $\Sel_{l^m} E_L \cong H^1(S^L, \bZ/l^m\bZ \oplus \mu_{l^m})$ admits a map to $\Hom(\Pic(S^L)/l^m \Pic(S^L), \bZ/l^m\bZ) \oplus \Pic(S^L)[l^m] $ with kernel of bounded size. 
\epf

\brems
\remi 
We expect that the conclusion of \Cref{S-T} holds already with $n(l, m) = 1$. 

\remi
For a composite $l_1^{m_1}\cdot \dotsc\cdot l_k^{m_k}$ prime to $p$, the proof gives a $q = p^{n(l_1, m_1, \dotsc, l_k, m_k)}$ such that for every finite extension $K/\bF_q(t)$, the unbounded growth of the number of $\bZ/l_i^{m_i}\bZ$-summands of $\Pic(S^L)[l_i^{m_i}]$ is simultaneous as $L/K$ ranges over quadratic extensions (of the form $L = L\pr K$ as in \Cref{S-T}).

\remi
A possible choice for $n(l, m)$ is $2n$ with $(-p)^n \equiv 1 \bmod l^m$ (e.g., $n(l, m) \ce 2l^{m - 1}(l - 1)$): in the proof take the supersingular $E \ra \Spec \bF_p$ which has $x^2 + p$ as the characteristic polynomial of the $p$-power Frobenius $\Frob_p$, so $\Frob_{p^{2n}}$ fixes $E[l^m]$.
\erems

\bcor \lab{sel-unb-ff}
If $l \neq p$, then there is a finite extension $K\pr/K$ (depending on $l$ and $A$) such that the analogue of \Cref{conj-sel} holds for $A_{K^{\prime\prime}}$ and $l$ for every finite extension $K^{\prime\prime}/K\pr$, i.e., $\#\Sel_l A_L$ is unbounded as $L/K^{\prime\prime}$ ranges over quadratic extensions.
\ecor

\bpf
Due to Theorems \ref{bdd-deg}\ref{bdd-deg-a} and \ref{S-T}, it suffices to choose $K\pr$ to contain $\bF_{p^n}$ with $n = n(l, 1)$ and satisfy either $\bZ/l\bZ \subset A[l]_{K\pr}$ or $\mu_l \subset A[l]_{K\pr}$.
\epf

\section{Iwasawa theory of class groups and Selmer groups} \lab{Iw-prel}

To keep the discussion focused, we assume in this and the next two sections that $K$ is a number field, even though the question of function field analogues is an interesting one. Likewise, we set aside the possibility of more general $p$-adic Lie extensions and fix a $\bZ_p$-extension $K_\infty/K$. Concretely, $K_\infty/K$ is Galois with $\Gal(K_\infty/K)\cong \bZ_p$; we fix a choice of the latter isomorphism, which identifies the Iwasawa algebra $\Lambda$ of $K_\infty/K$ with $\bZ_p[[T]]$. We denote by  $v_1, \dotsc, v_k$ the places of $K$ ramified in $K_\infty$, so $k \ge 1$ and $v_i \mid p$, and by $K_n$ the subfield of $K_\infty$ fixed by $p^n\bZ_p$.

\bpp[Iwasawa theory of class groups] \lab{Iw-cl}
Let $M$ be the maximal unramified abelian pro-$p$ extension of $K_\infty$. Set $X \ce \Gal(M/K_\infty)$, which is a finitely generated torsion $\Lambda$-module (cf.~\cite{Ser58}*{Thm.~5 et \S5}). The structure theory of such $\Lambda$-modules gives a $\Lambda$-homomorphism
\[
 X \ra \bigoplus_{i} \Lambda/f_i^{l_i}\Lambda \oplus \bigoplus_{j} \Lambda/p^{m_j}\Lambda,
\]
with finite kernel and cokernel (i.e., a \emph{pseudo-isomorphism}) for uniquely determined $m_j \in \bZ_{> 0}$, monic polynomials $f_i \in \bZ_p[[T]]$ that are monomials mod $p$, and $l_i \in \bZ_{> 0}$. The $\lambda$- and $\mu$-invariants of $K_\infty/K$ are 
\[ \lambda_{\Pic} \ce \sum l_i\deg f_i, \quad \quad \mu_{\Pic} \ce \sum m_j.\]
We also set $\mu_{\Pic}^{(m)} \ce \sum_{j} \min(m_j, m)$ for $m \ge 0$, which is of interest because it governs the growth of $\#\Pic(S^{K_n})[p^m]$ (as opposed to the customary in Iwasawa theory $\#\Pic(S^{K_n})[p\I]$): 
\epp

\bprop \lab{Iw-cl-gr}
$\#\Pic(S^{K_n})[p^m] \sim_{K, K_\infty, m}  p^{\mu_{\Pic}^{(m)}p^n}$.
\eprop

Before giving the proof we record a trivial lemma that clarifies implicit computations in subsequent arguments involving pseudo-isomorphisms; the lemma will be used without explicit notice.

\blem \lab{no-cite}
Let $R$ be a commutative ring, and let $X \xra{f} Y$ be a homomorphism of $R$-modules with finite kernel and cokernel. For $r \in R$, the induced $X/rX \xra{f_{/r}} Y/rY$ and $X[r] \xra{f_{[r]}} Y[r]$ satisfy
\[ \ba
\#\Ker f_{/r} &\le \#\Ker f \cdot \#\Coker f, \quad \quad \#\Coker f_{/r} \le \#\Coker f, \\
\#\Ker f_{[r]} &\le \#\Ker f, \quad \quad \quad \quad \quad \quad \quad  \#\Coker f_{[r]} \le \#\Ker f \cdot \#\Coker f.
\ea \]
\elem

\bpf 
Apply the snake lemma twice.
\epf

\bpf[Proof of \Cref{Iw-cl-gr}]
Replacing $K$ by $K_n$ has the effect of multiplying $\mu_{\Pic}^{(m)}$ by $p^n$ (since $\bZ_p[[T]]$ is replaced by $\bZ_p[[(T + 1)^{p^n} - 1]]$). By choosing $n$ large, we are therefore reduced to the case when each $v_i$ is totally ramified in $K_\infty$. 

In this case, by \cite{Ser58}*{Thm.~4}, as $\bZ_p$-modules, $\Pic(S^{K_n})[p\I]$ is isomorphic to the quotient of the finitely generated $X/((T + 1)^{p^n} - 1)X$ by a submodule generated by $k$ elements. Hence 
\[ 
\#\Pic(S^{K_n})[p^m] \sim \#(X/((T + 1)^{p^n} - 1)X)[p^m] \sim \prod_j \#(\Lambda/(p^{m_j}, (T + 1)^{p^n} - 1))[p^m] = p^{\mu_{\Pic}^{(m)}p^n}. \qedhere
\] 
\epf

\bpp [Iwasawa theory of Selmer groups] \lab{Iw-Sel}
The \emph{$p\I$-Selmer group} of $A_{K_n}$ is 
\[
\Sel_{p\I} A_{K_n} \ce \varinjlim_m \Sel_{p^m} A_{K_n},
\]
and that of $A_{K_\infty}$ is 
\[
\Sel_{p\I} A_{K_{\infty}} \ce \varinjlim_n \Sel_{p\I} A_{K_n}.
\]
For the compact Pontryagin dual $X\pr \ce \Hom(\Sel_{p\I} A_{K_{\infty}}, \bQ_p/\bZ_p)$, one knows \epp

\bcl \lab{sel-fg}
The $\Lambda$-module $X\pr$ is finitely generated.
\ecl

\bpf
Fix a nonempty open $U \subset S[\f{1}{p}]$ for which $\cA_U \ra U$ is an abelian scheme. Finiteness of $H^1_\et(U, \cA[p])$ \cite{Mil06}*{II.2.13} implies that of $H^1_\et(U, \cA[p\I])[p]$: the exact sequences 
\[\ba
&0 \ra \cA[p]_U \ra \cA[p^{n} ]_U \xra{p} \cA[p^{n - 1}]_U \ra 0, \\
&0 \ra \cA[p^{n - 1}]_U \ra \cA[p^{n} ]_U \xra{p^{n - 1}} \cA[p]_U \ra 0
\ea\]
give $\#H^1_\et(U, \cA[p^n])[p] \le \#H^1_\et(U, \cA[p]) \cdot \#A(K)[p]$. Consequently, $H^1_\et(U, \cA[p\I])$ is $\bZ_p$-cofinitely generated. 

Let $U_\infty \ce \varprojlim U_{S^{K_n}}$ be the normalization of $U$ in $K_\infty$. Since $U_\infty/U$ is pro-(finite \'{e}tale Galois), the Hochschild--Serre spectral sequence
\[
H^i(\Gal(K_\infty/K), H^j_\et(U_{\infty}, \cA[p\I])) \implies H^{i + j}_\et(U, \cA[p\I])
\]
shows that $H^1_\et(U_{\infty}, \cA[p\I])^{\Gal(K_\infty/K)}$ is $\bZ_p$-cofinitely generated. Therefore, so is 
\[
(\Sel_{p\I} A_{K_\infty})^{\Gal(K_\infty/K)} \subset H^1_\et(U_{\infty}, \cA[p\I])^{\Gal(K_\infty/K)}.
\]
Pontryagin duality then gives the finiteness of $X\pr/(T, p)$, and it remains to invoke the relevant version of Nakayama's lemma \cite{Ser58}*{Lemme~4}.
\epf

\Cref{sel-fg} and the structure theory of finitely generated $\Lambda$-modules give a pseudo-isomorphism
\be \lab{p-i-sel}
X\pr \ra \Lambda^{\rho} \oplus \bigoplus_{s} \Lambda/f_s^{l_s\pr}\Lambda \oplus \bigoplus_{t} \Lambda/p^{m_t\pr}\Lambda
\ee
as in \ref{Iw-cl} (with similar uniqueness claims). However, unlike $X$, the $\Lambda$-module $X\pr$ need not be torsion, i.e., $\rho > 0$ is possible. As for class groups, set $\mu_{\Sel}^{(m)} \ce \sum_{t} \min(m_t\pr, m)$ for $m \ge 0$.

\bpp[Controlled growth] \lab{c-gr}
We say that \emph{the control theorem holds} for $A$ and $K_\infty$, if 
\[
\Sel_{p\I} A_{K_n} \ra (\Sel_{p\I} A_{K_\infty})^{\Gal(K_\infty/K_n)}\ \ \text{ for }n\ge 0
\]
has finite kernel and cokernel of order bounded independently of $n$. The first result of this type is due to Mazur \cite{Maz72}*{6.4 (i)}; it has subsequently been generalized by Greenberg \cite{Gre03}*{5.1}: potential good ordinary reduction of $A$ at all $v \mid p$ is sufficient for the control theorem to hold. Such results play a purely axiomatic role in our computations:
\epp

\bprop \lab{Iw-sel-gr}
$\#\Sel_{p^m} A_{K_n} \sim_{A, K_\infty, m} p^{(\rho m + \mu_{\Sel}^{(m)})p^n}$, if the control theorem holds for $A$ and $K_\infty$.
\eprop

To replace $\Sel_{p^m} A_{K_n}$ by $(\Sel_{p\I} A_{K_n})[p^m]$ we will need a quantitative version of \cite{BKLPR13}*{5.9}:

\blem \lab{sel-tors}
Let $A \ra \Spec K$ be an abelian variety over a global field, $p$ a prime, and $a, b \in \bZ_{> 0}$.
\benum
\item \lab{sel-tors-b} The kernel and cokernel of $\Sel_a A \ra (\Sel_{ab}A)[a]$ are of size at most $\#A[a](K)$.

\item \lab{sel-tors-c} The kernel and cokernel of $\Sel_{p^m} A \ra (\Sel_{p\I} A)[p^m]$ are of size at most $\#A[p^m](K)$.
\eenum
\elem

\bpf 
Part \ref{sel-tors-c} is obtained from \ref{sel-tors-b} by taking direct limits. As for \ref{sel-tors-b}, the cohomology sequence of $0 \ra A[a] \ra A[ab] \xra{a} A[b] \ra 0$ gives the kernel claim since $\f{\#A[a](K) \cdot \#A[b](K) }{\#A[ab](K)} \le \# A[a](K)$. Selmer groups consist of $H^1$-classes that vanish in every $H^1(K_v, A)$, so $\f{(\Sel_{ab}A)[a]}{\im(\Sel_a A)} \hra \f{H^1(K, A[ab])[a]}{\im(H^1(K, A[a]))}$, and the cokernel claim results from the injection $\f{H^1(K, A[ab])[a]}{\im(H^1(K, A[a]))} \hra \Ker\p{H^1(K, A[b]) \ra H^1(K, A[ab])}.$
\epf

\bpf[Proof of \Cref{Iw-sel-gr}]
By \Cref{sel-tors}\ref{sel-tors-c}, the control theorem, and Pontryagin duality,
\[
\#\Sel_{p^m} A_{K_n} \sim \# (\Sel_{p\I} A_{K_n})[p^m]\sim  \#(\Sel_{p\I}A_{K_{\infty}})^{\Gal(K_\infty/K_n)}[p^m] \sim \#\p{X\pr/(p^m, (T + 1)^{p^n} - 1)}.
\] 
Therefore, the desired conclusion results from \eqref{p-i-sel} (and \Cref{no-cite}).
\epf

\section{Relations between the Iwasawa invariants of Selmer groups and class groups} \lab{Iw}

We keep the setup of \S\ref{Iw-prel} and denote by $\ord_p$ the $p$-adic valuation normalized by $\ord_p p = 1$.

\bprop \lab{Zp-low}
Suppose that the control theorem holds for $A$ and $K_\infty$, and let $\Sigma$ be the set of finite places of $K$ that decompose completely in $K_\infty$.
\benum
\item \lab{Zp-low-a}
If $A[p^m]$ has $\bigoplus_i \bZ/p^{a_i}\bZ \oplus \bigoplus_j \mu_{p^{b_j}}$ with $a_i, b_j \ge 1$ as a $K$-subgroup, $p \neq 2$, and each $v\mid p$ is finitely decomposed in $K_\infty$, then
\[
 \rho m + \mu_{\Sel}^{(m)} \ge \sum_i \mu_{\Pic}^{(a_i)} + \sum_j \mu_{\Pic}^{(b_j)} - r_2\sum_j b_j - \sum_{v\in \Sigma} \ord_p\p{ \f{\#\Phi_v(\bF_v)}{\#(p^m\Phi_v)(\bF_v)}}.
\]

\item \lab{Zp-low-b}
If $A[p^m]$ has $\bigoplus_i \bZ/p^{a_i}\bZ$ with $a_i \ge 1$ as a $K$-subgroup and $A$ has semiabelian reduction at all $v \mid p$, then
\[
\rho m + \mu_{\Sel}^{(m)} \ge \sum_i \mu_{\Pic}^{(a_i)} - \sum_{v\in \Sigma} \ord_p\p{ \f{\#\Phi_v(\bF_v)}{\#(p^m\Phi_v)(\bF_v)}}.
\]
\eenum
\eprop

\bpf
We begin with some preliminary observations.
\begin{enumerate}[label={(\arabic*)}]
\item \lab{a-pic}
$\#\Pic (S^{K_n}[\f{1}{p}])[p^{a_i}] \sim \#\Pic (S^{K_n})[p^{a_i}]$ in \ref{Zp-low-a}, since $\#(S^{K_n}\setminus S^{K_n}[\f{1}{p}])$ is bounded.

\item \lab{a-cx}
The number of complex places of $K_n$ is $r_2p^n$.

\item \lab{a-tam}
Since $S^{K_n}[\f{1}{p}] \ra S[\f{1}{p}]$ is \'{e}tale, $\prod_{\substack{w\nmid p\infty \\ w\text{ not above } \Sigma}} \#\Phi_w^{K_n} \sim 1$ where $w$ denotes a place of $K_n$.

\item \lab{b-tam}
For a place $w$ of semiabelian reduction for $A_{K_n}$, one has $\f{\#\Phi_w^{K_n}(\bF_w)}{\#(p^m\Phi_w^{K_n})(\bF_w)} \le  \#\Phi_w^{K_n}[p^m] \le~p^{2mg}$ where the last step uses surjectivity of multiplication by $p^m$ on $(\cA^{K_n})^0(\ov{\bF}_w)$ and the consideration of the finite part \cite{EGAIV4}*{18.5.11 c)} of the quasi-finite separated $(\cA^{K_n}[p^m])_{\cO_w}$.
\eenum

Combining \Cref{Iw-sel-gr,Iw-cl-gr} with \Cref{low-bds} and using \ref{a-pic}--\ref{b-tam}, we get
\[\ba
 p^{(\rho m + \mu_{\Sel}^{(m)})p^n} &\gtrsim_{A, K_\infty, m} p^{\p{\sum_i \mu_{\Pic}^{(a_i)} + \sum_j \mu_{\Pic}^{(b_j)} - r_2 \sum_j b_j}p^n} \cdot \p{\prod_{v\in \Sigma} \f{\#\Phi_v(\bF_v)}{\#(p^m\Phi_v)(\bF_v)}}^{-p^n}\quad \text{ and} \\
  p^{(\rho m + \mu_{\Sel}^{(m)})p^n} &\gtrsim_{A, K_\infty, m} p^{\p{\sum_i \mu_{\Pic}^{(a_i)}}p^n} \cdot \p{\prod_{v\in \Sigma} \f{\#\Phi_v(\bF_v)}{\#(p^m\Phi_v)(\bF_v)}}^{-p^n}
\ea\]
in cases \ref{Zp-low-a} and \ref{Zp-low-b}, respectively; the claimed inequalities follow by taking $n$ large enough.
\epf

\brem
The control theorem can hold in presence of completely decomposed places of bad reduction for $A$, see \cite{Gre03}*{5.1}.
\erem

\bprop \lab{Zp-up}
Set $r \ce r_1$ if $p = 2$, and $r \ce 0$ if $p \neq 2$. Suppose that the control theorem holds for $A$ and $K_\infty$, and every place $v$ above $p$ or of bad reduction for $A$ is finitely decomposed in $K_\infty$. If $A[p^m]$ has a filtration by $K$-subgroups with subquotients isomorphic to $\bZ/p^{a_i}\bZ$ or $ \mu_{p^{b_j}}$ with $a_i, b_j \ge 1$, then
\[
\rho m + \mu_{\Sel}^{(m)} \le 2mg[K : \bQ] - r_2\sum_j b_j + \sum_i (\mu_{\Pic}^{(a_i)} + r) + \sum_j \mu_{\Pic}^{(b_j)},
\]
and also
\[
 \rho m + \mu_{\Sel}^{(m)} \le mg[K : \bQ] + (r_1 + r_2)\sum_j b_j + \sum_i (\mu_{\Pic}^{(a_i)} + r) + \sum_j (\mu_{\Pic}^{(b_j)} + r) - \sum_{\substack{\text{real }v \\ p = 2}} \ord_2(\#\pi_0(A(K_v))).
\]
\eprop

\bpf 
We begin with some preliminary observations.
\begin{enumerate}[label={(\arabic*)}]
\item \lab{a-real}
Each infinite place of $K$ is completely decomposed in $K_\infty$.

\item \lab{a-pic+}
$\#(\Pic_+(S^{K_n})/p^{a_i}\Pic_+(S^{K_n})) \le 2^{rp^n} \cdot \#\Pic (S^{K_n})[p^{a_i}]$.

\item \lab{a-Upic}
$\#\Pic (U_{S^{K_n}})[p^{b_j}] \sim_{A, K_\infty, m} \#\Pic (S^{K_n})[p^{b_j}]$, since $\#(S^{K_n}\setminus U_{S^{K_n}})$ is bounded.
\eenum

Combining \Cref{Iw-sel-gr,Iw-cl-gr} with \Cref{up-bds} applied to $U_{S^{K_n}}$, where $U$ is the largest open subscheme of $S[\f{1}{p}]$ for which $\cA_{U} \ra U$ is an abelian scheme, and using \ref{a-real}--\ref{a-Upic}, we get
\[\ba
 p^{(\rho m + \mu_{\Sel}^{(m)})p^n} &\lesssim_{A, K_\infty, m} p^{\p{\sum_i (\mu_{\Pic}^{(a_i)} + r) + \sum_j \mu_{\Pic}^{(b_j)} + [K: \bQ]\sum_i a_i  + (r_1 + r_2)\sum_j b_j}p^n} \quad \text{ and} \\
  p^{(\rho m + \mu_{\Sel}^{(m)})p^n} &\lesssim_{A, K_\infty, m} p^{\p{\sum_i (\mu_{\Pic}^{(a_i)} + r) + \sum_j (\mu_{\Pic}^{(b_j)} + r) + mg[K : \bQ] + (r_1 + r_2)\sum_j b_j}p^n} \p{\prod_{\substack{\text{real }v \\ p = 2}} \#\pi_0(A(K_v))}^{-p^n}.
\ea\]
The claimed inequalities follow by taking $n$ large enough.
\epf

\bcor \lab{rho-bd}
Suppose that the control theorem holds for $A$ and $K_\infty$, and every place $v$ above $p$ or of bad reduction for $A$ is finitely decomposed in $K_\infty$. If $A[p]$ has a filtration by $K$-subgroups with $a$ subquotients isomorphic to $\bZ/p\bZ$ and $b$ subquotients isomorphic to $\mu_p$ with $a + b = 2g$, then
\[
\rho \le g[K : \bQ] + ar + 2g \mu_{\Pic}^{(1)} + \min\p{g[K : \bQ] - r_2b,\  b(r_1 + r_2 + r) }.
\]
\ecor

\bpf
Since $0 \subset A[p] \subset A[p^2] \subset \dotsc \subset A[p^m]$ has subquotients $A[p]$, \Cref{Zp-up} applies.
\epf

\brem
Remark \ref{ch-b} applies equally well to \Cref{Zp-up} and \Cref{rho-bd}.
\erem

\bpp \lab{K-prime}
The assumptions on $A[p^m]$ in \Cref{Zp-low,Zp-up} are satisfied after replacing $K$ by a finite extension $K\pr$. We record how this affects the Iwasawa invariants involved in the obtained inequalities. Set $K\pr_\infty \ce K\pr K_\infty$, and write $K_n\pr$, $\mu_{\Pic}\pr$, $\rho\pr$, $\mu_{\Sel}\pr$, etc.~for $K_\infty\pr/K\pr$ analogues of the familiar notation.
\epp

\blem \lab{mu-pic-gr}
One has $\mu_{\Pic}^{(m)} \le \mu_{\Pic}^{\prime (m)}$ for all $m\ge 0$. In particular, $\mu_{\Pic} \le \mu_{\Pic}\pr$.
\elem

\bpf
If $K\pr \cap K_\infty = K_{n}$, then \Cref{pic-extn} and \Cref{Iw-cl-gr} give $p^{n} \mu_{\Pic}^{(m)} \le \mu_{\Pic}^{\prime(m)}$.
\epf

\blem \lab{mu-sel-gr}
Suppose that the control theorem holds for $A$ and $K_\infty$ and also for $A_{K\pr}$ and $K\pr_\infty$. Then $\rho m + \mu_{\Sel}^{(m)} \le \rho\pr m + \mu_{\Sel}^{\prime(m)}$ for all $m\ge 0$. In particular, $\rho \le \rho\pr$, and if $\rho\pr = 0$, then $\mu_{\Sel} \le \mu_{\Sel}\pr$.
\elem

\bpf
If $K\pr \cap K_\infty = K_{n}$, then \Cref{sel-extn} and \Cref{Iw-sel-gr} give $p^{n}(\rho m + \mu_{\Sel}^{(m)}) \le \rho\pr m + \mu_{\Sel}^{\prime(m)}$.
\epf

\section{Conclusions for the cyclotomic $\bZ_p$-extension}\lab{cycl-concl}

Keeping the setup of \S\ref{Iw-prel}, we now assume that $K_\infty/K$ is the cyclotomic $\bZ_p$-extension, i.e., the unique $\bZ_p$-subextension of $K(\mu_{p\I})/K$. No anomalies occur: every finite $v$ is finitely decomposed in $K_\infty$, and for a finite extension $K\pr/K$, the compositum $K_\infty\pr \ce K\pr K_\infty$ is the cyclotomic $\bZ_p$-extension~of~$K\pr$.

\bconj[Iwasawa \cite{Iwa71}*{p.~392}, \cite{Iwa73}*{p.~11}] \lab{conj-Iw}
$\mu_{\Pic} = 0$. 
\econj

\bconj[Mazur \cite{Maz72}*{p.~184}] \lab{conj-Maz}
If $A$ has good ordinary reduction at all $v\mid p$, then $\rho = 0$.
\econj

\bpp[Status of \ref{conj-Iw} and \ref{conj-Maz}] \lab{status}
\Cref{conj-Iw} is known for abelian $K/\bQ$ \cite{FW79}; \Cref{conj-Maz} is known for $A = E_K$, if $p$ is odd, $E \ra \Spec \bQ$ is an elliptic curve with good ordinary reduction at $p$, and $K/\bQ$ is abelian \cite{Kat04}*{17.4}, \cite{Roh84}, and also for $A$ with finite $\Sel_{p\I} A$, as the control theorem shows. Examples with $\mu_{\Sel} > 0$ are known, and in fact $\mu_{\Sel}^{(1)}$ can be arbitrarily large when $K$ is allowed to vary, as \Cref{large-mu1} shows.
\epp

The inequalities of \S\ref{Iw} allow one to relate \Cref{conj-Iw,conj-Maz}:

\bthm \lab{main-Iw}
If $\rho + \mu_{\Sel} = 0$, the control theorem holds for $A$ and $K_\infty$, and
\begin{enumerate}[label={(\roman*)}]
\item 
$A$ has $\bZ/p\bZ$ as a $K$-subgroup and semiabelian reduction at all $v\mid p$, or

\item $p$ is odd, and $A$ has $\bZ/p\bZ$ as a $K$-subgroup, or

\item 
$p$ is odd, $K$ is totally real, and $A$ has $\mu_p$ as a $K$-subgroup,
\eenum
then $\mu_{\Pic} = 0$.
\ethm

\bpf
The conclusion is immediate from \Cref{Zp-low}, because $\Sigma = \emptyset$.
\epf

Adopting the notation of \ref{K-prime}, one can use the results of \S\ref{Iw} to study boundedness questions:

\bthm
If $K = \bQ$, the reduction of $A$ at $p$ is good ordinary, and $A[p]$ has a filtration by $K$-subgroups with subquotients isomorphic to $\bZ/p\bZ$ or $\mu_p$, then $\rho\pr \lesssim_{d, g} 1$ and $\mu_{\Sel}^{\prime(1)} \lesssim_{d, g} 1$ for an abelian extension $K\pr/\bQ$ of degree $d$.
\ethm

\bpf
Indeed, $\mu_{\Pic}\pr = 0$ (cf.~\ref{status}), so \Cref{Zp-up} gives the claim.
\epf

\brems
\remi
If one assumes \Cref{conj-Iw}, then the abelian restriction on $K\pr/\bQ$ is not needed; in fact, one can then also drop the assumption on $A[p]$ and get the conclusion $\rho\pr, \mu_{\Sel}^{\prime(1)} \lesssim_{d, g, p} 1$ with the help of \Cref{mu-sel-gr}. Conversely, due to \Cref{Zp-low} and \Cref{mu-sel-gr}, such a conclusion for all $d$ and a single $A$ with good ordinary reduction at $p$  would give $\mu_{\Pic}^{\prime (1)} \lesssim_{d, p} 1$. Due to \Cref{Zp-up} and \Cref{mu-pic-gr}, this would in turn imply $\rho\pr, \mu_{\Sel}^{\prime(1)} \lesssim_{d, g, p} 1$ for \emph{every} $A$ with good ordinary reduction at $p$. Is there a way to prove $\rho\pr, \mu_{\Sel}^{\prime(1)} \lesssim_{d, g, p} 1$ for a single such $A$ without restricting to abelian $K\pr/\bQ$ and relying on \Cref{conj-Iw}?

\remi
If $d = g = 1$ and the reduction of $A$ at $p$ is good ordinary (but no assumption on $A[p]$), then Greenberg has conjectured that $\mu_{\Sel}^{\prime (1)} \le 1$ \cite{Gre99}*{1.11 and p.~118 Remark}. We show that $\mu_{\Sel}^{\prime (1)}$ can grow unboundedly as $d$ grows:
\erems

\beg \lab{large-mu1}
Suppose that $\cA[p] \cong (\bZ/p\bZ)^g \oplus \mu_p^g$ over $S$ and $A$ has good reduction at all $v\mid p$. Then $\cA^{K\pr}[p] \cong (\bZ/p\bZ)^g \oplus \mu_p^g$ over $S^{K\pr}$ for every finite extension $K\pr/K$ \cite{Ces13}*{3.4 and the proof~of~3.3}. For instance, this is the case for $K = \bQ$ and $A = X_0(11)$ with $p = 5$ \cite{Ces13}*{1.12}.

Assume that $p > 2$. By \cite{Ces13}*{5.5 and the proof of 5.4} and \Cref{const,mu},
\be \lab{sel-comp}
\#\Sel_p A_{K_n} \sim_{A, K_\infty} \#H^1(S^{K_n}, (\bZ/p\bZ)^g \oplus \mu_p^g) \sim \#\Pic(S^{K_n})[p]^{2g} \cdot p^{gp^n(r_1 + r_2)}.
\ee
If the reduction is ordinary at all $v\mid p$, then \eqref{sel-comp} combines with \Cref{Iw-cl-gr,Iw-sel-gr} to give
\[
\rho + \mu_{\Sel}^{(1)} = 2g\mu_{\Pic}^{(1)} + g(r_1 + r_2).
\]
The same reasoning applies with $K$ replaced by a finite extension $K\pr$. In particular, if $p > 2$, the reduction of $A$ at all $v\mid p$ is good ordinary, and $\cA[p] \cong (\bZ/p\bZ)^g \oplus \mu_p^g$, then
\[
\rho\pr + \mu_{\Sel}^{\prime (1)} = 2g\mu_{\Pic}^{\prime (1)} + g(r_1\pr + r_2\pr)
\]
for every finite extension $K\pr/K$. In particular, under \Cref{conj-Iw,conj-Maz}, $\mu_{\Sel}^{\prime (1)} = g(r_1\pr + r_2\pr)$, and for $K = \bQ$, $A = X_0(11)$, $p = 5$, and $K\pr/\bQ$ abelian, the same holds unconditionally (cf.~\ref{status}).
\eeg

\appendix

\section{Cardinalities of the images of local Kummer homomorphisms} \lab{app-a}

Let $K$ be a local field, $A$ a $g$-dimensional abelian variety over $K$, and $l$ a prime. \Cref{card-kum} summarizes standard computations in the form needed for the bounds of \S\S\ref{low}--\ref{up}.

\bprop\lab{card-kum} 
Fix an $m\in \bZ_{> 0}$. If $K$ is nonarchimedean, let $\bF_K$ be its residue field.
\benum
\item \lab{card-kum-a}
If $K$ is nonarchimedean and $l \neq \Char \bF_K$, then $\#\p{A(K)/l^mA(K)} = \#A(K)[l^m]$.

\item \lab{card-kum-b}
If $K$ is a finite extension of $\bQ_l$, then $\#\p{A(K)/l^mA(K)} = l^{mg[K: \bQ_l]}\cdot \#A(K)[l^m]$.

\item \lab{card-kum-c} 
If $K \cong \bR$ and $l = 2$, then $A(K)/l^mA(K) \cong \pi_0(A(K))$ (component group for the archimedean topology) and $\#\pi_0(A(K)) \le 2^{g}$. In all other archimedean cases, $A(K)/l^mA(K) = 0$.

\eenum
\eprop

\bpf \hfill
\benum
\item 
Let $\cO_K$ be the ring of integers of $K$ and $\cA \ra \Spec \cO_K$ the N\'{e}ron model of $A$. Since $\cA$ is smooth over the Henselian $\cO_K$, the reduction homomorphism $\cA(\cO_K) \ra \cA(\bF_K)$ is surjective \cite{BLR90}*{2.2/14}; once we show that its kernel is uniquely divisible by $l^m$, the conclusion follows from the snake lemma because $\#(\cA(\bF_K)/l^m\cA(\bF_K)) = \#\cA(\bF_K)[l^m]$ due to finiteness of $\cA(\bF_K)$. Since $\cA \xra{l^m} \cA$ is separated \'{e}tale \cite{BLR90}*{7.3/2(b)}, so is its pullback over each $P \in \cA(\cO_K)$, and the claimed unique divisibility follows from \cite{EGAIV4}*{18.5.12}.

\item
The finite index inclusion $\bZ_l^{g[K : \bQ_l]} \subset A(K)$ of \cite{Mat55}*{Thm.~7} with the snake lemma give
\[ 
\f{\#\Coker\p{A(K) \xra{l^m} A(K)}}{\# \Ker\p{A(K) \xra{l^m} A(K)}} = \f{\#\Coker\p{\bZ_l^{g[K: \bQ_l]} \xra{l^m} \bZ_l^{g[K: \bQ_l]}}}{\# \Ker\p{\bZ_l^{g[K: \bQ_l]} \xra{l^m} \bZ_l^{g[K: \bQ_l]}}} = l^{mg[K : \bQ_l]}. 
\]

\item 
$H^1(K, A[l^m]) = 0$ unless $K \cong \bR$ and $l = 2$, in which case \cite{GH81}*{1.1 (3)} applies. \qedhere
\eenum
\epf

\brem \lab{pos-inf}
Finiteness of quotients $A(K)/l^mA(K)$ fails for $K$ of characteristic $l$: for instance, for the Tate elliptic curve $\bG_m/q^{\bZ}$, combine the snake lemma with the well-known infinitude of $K^{\times} /K^{\times l^m}$ \cite{Iwa86}*{(2.2) and 2.8}.
\erem

\section{The flat cohomology of $\bZ/l^a\bZ$ and $\mu_{l^b}$} \lab{app-b}

Fix a nonempty open $U \subset S$. We work out the cardinalities of the (compactly supported) flat cohomology groups of $U$ with $\bZ/l^a\bZ$ or $\mu_{l^b}$ coefficients, which are needed in \S\S\ref{low}--\ref{up}

\blem \lab{const} \hfill
\benum
\item \lab{const-a} $\#H^1(S, \bZ/l^a\bZ) = \#(\Pic_+ S / l^a \Pic_+ S)$.

\item \lab{const-b}
If $l \neq \Char K$, then, interpreting $[K_v : \bQ_l]$ as $0$ unless $\Char K = 0$ and $v\mid l$,
 \[
 \#(\Pic_+ S / l^a \Pic_+ S) \le  \#H^1(U, \bZ/l^a\bZ) \le \#(\Pic_+ S / l^a \Pic_+ S) \cdot \prod_{v\in S \setminus U}\p{ \#\mu_{l^a}(K_v) \cdot l^{a[K_v: \bQ_l]}}.
 \]
\eenum
\elem

\bpf \hfill
\benum
\item Since $H^1(S, \bZ/l^a\bZ) \cong \Hom(\pi_1^\et(S), \bZ/l^a\bZ)$, 
the theory of the narrow Hilbert class field gives the claim in the number field case. For function fields, one can (alternatively) use duality: by \cite{Mil06}*{III.8.2}, $H^1(S, \bZ/l^a\bZ) \cong H^2(S, \mu_{l^a})^*$, so, due to $0 \ra \mu_{l^a} \ra \bG_m \xra{l^a} \bG_m \ra 0$ that is exact in $S_\fppf$, the vanishing of the Brauer group of $S$ gives the claim.

\item The exact $0 \ra H^1(S, \bZ/l^a\bZ) \ra H^1(U, \bZ/l^a\bZ) \ra \prod_{v \in S \setminus U} H^1(K_v, \bZ/l^a\bZ)/H^1(\cO_v, \bZ/l^a\bZ)$, \ref{const-a}, and local class field theory give the bounds, because $H^1(K_v, \bZ/l^a\bZ) \cong \Hom(K_v^\times/K_v^{\times l^a}, \bZ/l^a\bZ)$ and $H^1(\cO_v, \bZ/l^a\bZ) \cong \Hom(\pi_1^\et(\cO_v), \bZ/l^a\bZ) \cong \bZ/l^a\bZ$.    \qedhere
\eenum
\epf

\blem \lab{mu} 
$\#H^1(U, \mu_{l^b}) = \#\Pic(U)[l^b] \cdot l^{b \cdot \max(r_1 + r_2 + \#(S\setminus U) - 1, 0)} \cdot \#\mu_{l^b}(K)$. 
\elem

\bpf
Since $0 \ra \mu_{l^b} \ra \bG_m \xra{l^b} \bG_m \ra 0$ is exact in $U_\fppf$, its long exact cohomology sequence together with the unit theorem \cite{AW45}*{p.~491, Thm.~6} give the claim.
\epf

\blem\lab{const-c} 
Set $r \ce r_1$ if $l = 2$, and $r \ce 0$ if $l \neq 2$.
\benum
\item \lab{const-c-a}
If $a \ge 1$, then $\# H^1_c(S, \bZ/l^a\bZ) = \#(\Pic S/l^a \Pic S) \cdot 2^{\max(r - 1, 0)}.$

\item \lab{const-c-b}
If $a \ge 1$ and $U \neq S$, then $\# H^1_c(U, \bZ/l^a\bZ) = \#(\Pic U/l^a \Pic U) \cdot l^{a(\#(S \setminus U) - 1)} \cdot 2^r.$
\eenum
\elem

\bpf 
By duality \cite{Mil06}*{III.3.2, III.8.2}, $\#H^1_c(U, \bZ/l^a\bZ) = \#H^2(U, \mu_{l^a})$, and the claim follows from the cohomology sequence of $0 \ra \mu_{l^a} \ra \bG_m \xra{l^a} \bG_m \ra 0$ since the Brauer group of $U$ is understood from the exact sequence $0 \ra \Br U \ra \bigoplus_{v\not\in U} \Br(K_v) \xra{\sum \inv_v} \bQ/\bZ$ \cite{Mil06}*{II.2.1}.
\epf

\blem\lab{mu-c} 
If $b \ge 1$ and $l$ is invertible on $U \neq S$, then
\[
\ba
\#H^1_c(U, \mu_{l^b}) &\ge \#(\Pic_+ S / l^b \Pic_+ S) \cdot l^{-b(r_2 + 1)},   \\
\#H^1_c(U, \mu_{l^b}) &\le \#(\Pic_+ S / l^b \Pic_+ S) \cdot l^{b(r_1 + r_2 - 1)} \cdot  \prod_{v\in S \setminus U} \#\mu_{l^b}(K_v).
\ea
\]
 \elem
 
 \bpf 
We replace compactly supported flat cohomology by its \'{e}tale counterpart \cite{Mil06}*{II.2.3 and the preceding subsection}: by \cite{Mil06}*{II.3.3, III.3.2, III.8.1} and \cite{Gro68}*{11.7 1$^{\circ}$)}, the two meanings of $H^i_c(U, \mu_{l^b})$ agree.
 
By the Euler characteristic formula \cite{Mil06}*{II.2.13 (b)} and duality \cite{Mil06}*{II.3.3},
 \[ 
 \#H^1_c(U, \mu_{l^b}) = \f{\#H^0_c(U, \mu_{l^b}) \cdot \#H^2_c(U, \mu_{l^b}) }{2^{r}l^{br_2}\cdot  \#H^3_c(U, \mu_{l^b})} = \f{\#H^0_c(U, \mu_{l^b}) \cdot \#H^1(U, \bZ/l^b\bZ) }{ 2^{r}l^{b(r_2 + 1)}} 
 \] 
with $r$ as in \Cref{const-c}. By \cite{Mil06}*{II.2.3 (a)} (we use the $U \neq S$ assumption to discard $H^0(U, \mu_{l^b})$), 
\[
H^0_c(U, \mu_{l^b}) \cong \bigoplus_{v\mid \infty} \wh{H}^{-1}(K_v, \mu_{l^b}) \cong \bigoplus_{\text{real } v} H^{1}(K_v, \mu_{l^b}) \cong \bigoplus_{\text{real } v} \bR^{\times}/\bR^{\times l^b}
\]
where $\wh{H}^i$ denotes Tate cohomology. It remains to take into account \Cref{const}\ref{const-b}.
\epf


\end{document}